\documentclass[11pt]{amsart}
\usepackage{amsfonts}
\usepackage{amssymb}
\usepackage{amsmath}  
\usepackage{verbatim}
\usepackage{stmaryrd}

\input xy
\xyoption{all}

\theoremstyle{plain} 
\newtheorem{theorem}{Theorem}[section]
\newtheorem{lemma}[theorem]{Lemma}

\newtheorem{proposition}[theorem]{Proposition}
\newtheorem{corollary}[theorem]{Corollary}

\newtheorem{conjecture}[theorem]{Conjecture}

\theoremstyle{definition} 

\newtheorem{definition}[theorem]{Definition}

\newcommand{\End}[1]{\operatorname{\rm End}_{#1}}

\newcommand{\Hom}[1]{\operatorname{{\rm Hom}}_{#1}}

\newcommand{\dimv}{\underline{\dim}\,}

\newcommand{\add}{\mbox{{\rm add \!}}}

\newcommand{\MOD}{\mbox{{\rm mod \!}}}

\newcommand{\ind}[1]{\operatorname{\rm ind}_{#1}}

\newcommand{\Gr}[1]{\mbox{{\rm Gr}}_{#1}}

\newcommand{\demo}[1]{\textsc{Proof.} #1 \hfill $\Box$ \bigskip}

\newcommand{\cA}{\mathcal{A}}

\newcommand{\cC}{\mathcal{C}}
\newcommand{\cD}{\mathcal{D}}

\newcommand{\cF}{\mathcal{F}}

\newcommand{\bP}{\mathbb{P}}
\newcommand{\bQ}{\mathbb{Q}}
\newcommand{\bZ}{\mathbb{Z}}

\newcommand{\be}{\mathbf{e}}

\newcommand{\bi}{\mathbf{i}}

\newcommand{\bu}{\mathbf{u}}
\newcommand{\bx}{\mathbf{x}}
\newcommand{\by}{\mathbf{y}}

\begin{document}

\title[Linear independence of cluster monomials]{Linear independence of cluster monomials for skew-symmetric cluster algebras}
\author[G. Cerulli \and B. Keller \and D. Labardini \and P.-G. Plamondon]{Giovanni Cerulli Irelli \and Bernhard Keller \and Daniel Labardini-Fragoso \and Pierre-Guy Plamondon}

\address{Mathematisches Institut \\
         Universit\"at Bonn \\
          53115 Bonn \\
           Germany}
\address{ Universit\'e Paris Diderot - Paris 7 \\
          UFR de Math\'ematiques, Case 7012 \\
          B\^atiment Chevaleret \\
          75205 Paris Cedex 13 \\
          France}
\address{ Mathematisches Institut \\
         Universit\"at Bonn \\
          53115 Bonn \\
           Germany}
\address{ Laboratoire LMNO \\
          \'Equipe Alg\`ebre, G\'eom\'etrie et Logique \\
          Universit\'e de Caen \\
          F14032 Caen Cedex \\
          France}

\email{cerulli@math.uni-bonn.de}
\email{keller@math.jussieu.fr}
\email{labardini@math.uni-bonn.de}
\email{pierre-guy.plamondon@unicaen.fr}

\thanks{The fourth author was financially supported by the FQRNT}


\begin{abstract}
Fomin-Zelevinsky conjectured that in any cluster algebra,
the cluster monomials are linearly independent and that
the exchange graph and cluster complex are independent of the choice of
coefficients. We confirm these conjectures for all
skew-symmetric cluster algebras.
\end{abstract}

\dedicatory{Dedicated to Idun Reiten on the occasion of her seventieth birthday.}

\maketitle

\tableofcontents

\section{Introduction}

When S.~Fomin and A.~Zelevinsky invented cluster algebras
more than ten years ago \cite{FominZelevinsky02}, one of their
hopes was that the cluster monomials would be part of a
`canonical' basis. This naturally led them to the conjecture that 
cluster monomials should be linearly independent 
(Conjecture~4.16 of \cite{FominZelevinsky03a}).
The conjecture has been checked for several classes, notably
for finite type cluster algebras \cite{CK08} \cite{Demonet11},
for skew-symmetric acyclic cluster algebras 
\cite{GeissLeclercSchroeer10a} and for cluster algebras
associated with surfaces \cite{CL11} (see also \cite{FominShapiroThurston08} \cite{MSW11}). It was proved for all
skew-symmetric cluster algebras of geometric type whose
exchange matrix is of full rank in \cite{DWZ09} and \cite{Plamondon10}.
In this note, we show that the condition on the rank may
be dropped, and that we need not restrict ourselves to the 
geometric type. In other words, we show that cluster monomials
are linearly independent in any cluster algebra associated
with a skew-symmetric matrix. 

Our proof is based on the additive categorification of cluster
algebras using the cluster category \cite{Plamondon10}
(closely related to that using decorated representations
\cite{DWZ09}) and on an idea from \cite{Cerulli11}. We show that
the canonical cluster character yields a bijection from
the set of isomorphism classes of reachable rigid objects
onto the set of cluster monomials (Corollary~\ref{coro::bijection}). In this way,
we treat the cluster algebras of geometric type, \emph{i.e.} those associated
with ice quivers. As a corollary 
of our proof, we show that for cluster algebras associated with 
ice quivers, the exchange graph (and, in fact, the cluster complex) does not depend on the frozen part
of the quiver, which confirms Conjecture 4.3 of \cite{FZ07} for
this class.  This has been shown previously in \cite{GSV08} (see also the book \cite{GSV10}) for skew-symmetrizable cluster algebras whose defining matrix is of full rank. Once these results have been obtained for 
cluster algebras of geometric type,
we show how we can extend them to arbitrary coefficients
thanks to Fomin-Zelevinsky's separation formula.

We expect that our methods and results can be extended to the
skew-symmetrizable cluster algebras categorified by
L.~Demonet in \cite{Demonet10}.

\section{Recollections}
\subsection{Skew-symmetric cluster algebras of geometric type}
We recall notions from \cite{FZ07}.  Let $\widetilde{B}$ be an $n\times r$ matrix whose top $r\times r$ matrix $B$ is skew-symmetric.  We will work over the field $\bQ(x_1, \ldots, x_n)$; let $\bP$ be the tropical semifield generated by $x_{r+1}, \ldots, x_n$.

\begin{definition}
A \emph{seed} is a pair $(\widetilde{B}, \bu)$, where
\begin{itemize}
	\item $\widetilde{B}$ is a matrix as above, and
	\item $\bu = (u_1, \ldots, u_n)$ is an ordered set whose elements form a free generating subset of the field $\bQ(x_1, \ldots, x_n)$.
\end{itemize}
\end{definition}

\begin{definition}
Let $i \in \{1, \ldots, r\}$.  The \emph{mutation} of the seed $(\widetilde{B}, \bu)$ at vertex $i$ is the seed $(\widetilde{B'}, \bu')$, where
\begin{itemize}
	\item $\widetilde{B'} = (b'_{kj})$ is the $n\times r$ matrix defined by
	   \begin{equation*}
b'_{kj} =
\begin{cases}
-b_{kj} & \text{if $k=i$ or $j=i$},\\
b_{kj} + \frac{b_{ki} |b_{ij}| + |b_{ki}| b_{ij}}{2} & \text{otherwise},
\end{cases}
\end{equation*}
where $[z]_+ = \max(0, z)$ for any real number $z$.
	
	\item $\bu' = (u_1', \ldots, u_n')$, with $u_j' = u_j$ if $i\neq j$ and
	  \begin{displaymath}
	    u_iu_i' = \prod_{j=1}^{n}u_j^{[b_{ji}]_+} + \prod_{j=1}^{n}u_j^{[-b_{ji}]_+}.
    \end{displaymath}
    
\end{itemize}
\end{definition}

\begin{definition}
Let $\widetilde{B}$ be a matrix as above.
\begin{enumerate}
	\item  The \emph{clusters} are the sets $\bu$ contained in seeds $(\widetilde{B'}, \bu)$ obtained from the \emph{initial seed} $(\widetilde{B}, \bx = (x_1, \ldots, x_n))$ by iterated mutations.
	\item The \emph{cluster variables} are elements $u_i$ of any cluster $\bu$, with $i\in \{1, \ldots, r\}$.
	\item The \emph{coefficients} are the elements $x_{r+1}, \ldots, x_n$.
	\item A \emph{cluster monomial} is a product of cluster variables all belonging to the same cluster.
	\item The \emph{(skew-symmetric) cluster algebra (of geometric type)} $\cA(\widetilde{B})$ is the $\bZ\bP$-subalgebra of $\bQ(x_1, \ldots, x_n)$ generated by all cluster variables.
\end{enumerate}
\end{definition}

Notice that one can associate an ice quiver $(Q,F)$ to the matrix $\widetilde{B}$ in the following way:
\begin{itemize}
	\item $Q$ has $n$ vertices, labelled by the rows of $\widetilde{B}$;
	\item the set $F$ of frozen vertices is the set $\{r+1, \ldots, n\}$;
	\item for any vertices $i$ and $j$ of $Q$, there are $[b_{ij}]_+$ arrows from $i$ to $j$.
\end{itemize}

\begin{conjecture}[Conjecture 7.2 of \cite{FZ07}]
The cluster monomials of any cluster algebra are linearly independent over $\bZ\bP$.
\end{conjecture}

\subsection{Cluster categories}
We now briefly recall results on the categorification of cluster algebras via triangulated categories.  Recent surveys on the subject include \cite{R08}, \cite{R10}, \cite{Amiot11} and \cite{KellerSurvey12}.

Let $(Q,W)$ be a quiver with a non-degenerate potential \cite{DWZ08}.  Generalizing the cluster category introduced in \cite{BMRRT06}, C.~Amiot \cite{Amiot08} has defined the cluster category $\cC_{Q,W}$ of $(Q,W)$; it is a triangulated category with a rigid object $\Gamma = \Gamma_1 \oplus \ldots \oplus \Gamma_n$ whose endomorphism algebra is isomorphic to the Jacobian algebra $J(Q,W)$.

As in \cite[Definition 3.9]{Plamondon09}, define $\cD$ to be the full subcategory of all objects $X$ of $\cC_{Q,W}$ such that
\begin{itemize}
	\item there exists a triangle $T_1^X \rightarrow T_0^X \rightarrow X \rightarrow \Sigma T_1^X$, with $T_0^X$ and $T_1^X$ in $\add\Gamma$;
	\item there exists a triangle $T^0_X \rightarrow T^1_X \rightarrow \Sigma X \rightarrow \Sigma T^0_X$, with $T^0_X$ and $T^1_X$ in $\add\Gamma$; and
	\item $\Hom{\cC}(\Gamma, \Sigma X)$ is finite-dimensional.
\end{itemize}

Define the \emph{index} \cite{DK08} of an object $X$ of $\cD$ as
\begin{displaymath}
	\ind{\Gamma} X = [T_0^X] - [T_1^X] \in K_0(\add\Gamma) \cong \bZ^n.
\end{displaymath}
Note that $K_0(\add\Gamma)$ is isomorphic to the free abelian group generated by the elements corresponding to the indecomposable direct summands of $\Gamma$, which we denote by $[\Gamma_1], \ldots, [\Gamma_n]$.

An object $X$ \emph{admits an index} if there exists a triangle $T_1^X \rightarrow T_0^X \rightarrow X \rightarrow \Sigma T_1^X$, with $T_0^X$ and $T_1^X$ in $\add\Gamma$ (so that the index of $X$ is defined).

The \emph{cluster character} \cite{CC06} \cite{Palu08} is the map $CC: Obj(\cD) \longrightarrow \bZ[x_1^{\pm 1}, \dots, x_n^{\pm 1}]$ defined by
\begin{eqnarray*}
	CC(X) &=& \bx^{\ind{\Gamma}X}\sum_{\be}\chi\Big( \Gr{\be}\big( \Hom{\cC}(\Gamma, \Sigma X)   \big) \Big)\bx^{B \cdot \be},
\end{eqnarray*}
where \begin{itemize}
	       \item $B$ is the skew-symmetric matrix associated to $Q$;
	       \item $\Gr{\be}\big( \Hom{\cC}(\Gamma, \Sigma X)  \big)$ is the \emph{quiver Grassmannian} (see Section 2.3 of \cite{CC06}).  It is a projective variety whose points correspond to  submodules with dimension vector $\be$ of the $\End{\cC}(\Gamma)$-module $\Hom{\cC}(\Gamma, \Sigma X)$;
	       \item $\chi$ is the Euler--Poincar\'e characteristic.
      \end{itemize} 

A result of \cite{KR07} implies the following proposition if $H^0\Gamma$ is finite-dimensional.
\begin{proposition}[\cite{Plamondon09}]\label{prop::equivalence}
The functor $F = \Hom{\cC}(\Gamma, \Sigma(?))$ induces an equivalence of additive categories
\begin{displaymath}
	\xymatrix{\cD/(\Gamma) \ar[r]^-{\sim} & \MOD J(Q,W),}
\end{displaymath}
where $\MOD J(Q,W)$ is the category of finite-dimensional right modules over the Jacobian algebra $J(Q,W)$.
\end{proposition}

There is a notion of mutation of objects inside $\cC_{Q,W}$, studied in \cite{KY09}.  A \emph{reachable object} is an object obtained by iterated mutations from $\Gamma$.  Any reachable object lies in the subcategory $\cD$ \cite[Section 3.3]{Plamondon09}.

 The following theorem generalizes results obtained in \cite{CC06} \cite{CK08} \cite{Palu08}.
\begin{theorem}[\cite{Plamondon09}]\label{theo::surjection}
Let $B$ be an $r\times r$ skew-symmetric matrix, and let $Q$ be the associated quiver. Let $W$ be a non-degenerate potential on $Q$.  Then the cluster character $CC$ induces a surjection
\begin{displaymath}
	\{ \textrm{isoclasses of indecomposable reachable objects of $\cD$} \} \rightarrow \{ \textrm{cluster variables of $\cA(B)$} \},
\end{displaymath}
which commutes with mutation on both sides.
\end{theorem}

In particular, if $\widetilde{B}$ is an $n\times r$ matrix whose upper $r\times r$ matrix is $B$, and if $(\widetilde{Q},F)$ is the associated ice quiver, then the cluster variables of $\cA(\widetilde{B})$ are obtained from the cluster character $CC$ by restricting mutation on the left hand side to non-frozen vertices.

\section{Linear independence of cluster monomials}

\begin{definition}[Definition 6.1 of \cite{CL11}]
Let $B$ be a skew-symmetric matrix and $\widetilde{B}$ be an $n\times r$ integer matrix whose top $r\times r$ matrix is $B$.
\begin{enumerate}
	\item Let $\bu = (u_1, \dots, u_r)$ be a cluster of $\cA(\widetilde{B})$.  A \emph{proper Laurent monomial in $\bu$} is a product of the form $u_1^{c_1}\cdots u_r^{c_r}$ where at least one of the $c_i$ is negative.
	\item The cluster algebra $\cA(\widetilde{B})$ \emph{has the proper Laurent monomial property} if for any two clusters $\bu$ and $\bu'$ of $\cA(\widetilde{B})$, every monomial in $\bu'$ in which at least one element from $\bu' \setminus \bu$ appears with positive exponent is a $\bZ\bP$-linear combination of proper Laurent monomials in $\bu$.
\end{enumerate}
\end{definition}

\begin{theorem}[Theorem 6.4 of \cite{CL11}]
If $\cA(\widetilde{B})$ has the proper Laurent monomial property, then its cluster monomials are linearly independent over $\bZ\bP$.
\end{theorem}

Using the key idea of Lemma~5.1 in \cite{Cerulli11} it was proved in
\cite[Theorem 6.2]{CL11} that a cluster algebra $\cA(\widetilde{B})$ has the proper Laurent monomial property if there exists a potential on the quiver $Q$ of $\widetilde{B}$ satisfying certain conditions.  These conditions allow the use of a homological interpretation of the $E$-invariant (see \cite[Section 10]{DWZ09}).

If $(Q,W)$ is any quiver with a non-degenerate potential $W$, then the $E$-invariant has a similar interpretation in the cluster category $\cC_{Q,W}$ (see \cite[Proposition 4.15]{Plamondon10}): it can be thought of as one half of the dimension of the space of selfextensions of a given object.  Using it, we prove the following result.

Let $B$ be any skew-symmetric matrix and $\widetilde{B}$ be an $n\times r$ integer matrix whose top $r\times r$ matrix is $B$.  Let $(Q,F)$ be the ice quiver corresponding to $\widetilde{B}$, and let $W$ be a non-degenerate potential on $Q$.

\begin{theorem}\label{theo::rigid}
If $R$ is a rigid object of $\cD \subset \cC_{Q,W}$ which does not belong to $\add\Gamma$ and such that $\Hom{\cC}(\Gamma_i, \Sigma R) = 0$ for $i=r+1, \ldots, n$, then $CC(R)$ is a $\bZ\bP$-linear combination of proper Laurent monomials in the initial cluster.
\end{theorem}
\demo{ We follow the proof of \cite[Theorem 6.2]{CL11}.  Since $\{[\Gamma_1], \ldots, [\Gamma_n]\}$ is a $\bZ$-basis of $K_0(\add\Gamma)$, we can write $\ind{\Gamma}R = \sum_{i=1}^{n}g_i[\Gamma_i]$. Then we have that
\begin{eqnarray*}
	CC(R) &=& \bx^{\ind{\Gamma}R} \sum_{\be}\chi\Big( \Gr{\be} \big( \Hom{\cC}(\Gamma, \Sigma R) \big)  \Big) \bx^{\widetilde{B} \cdot \be_0} \\
	&=&  \sum_{\be}\chi\Big( \Gr{\be} \big( \Hom{\cC}(\Gamma, \Sigma R) \big)  \Big) \bx_N^{\ind{\Gamma}^N R + B \cdot \be_0} \bx_F^{\ind{\Gamma}^F R + B_c \cdot \be_0},
\end{eqnarray*}
where $B_c$ is the lower $(n-r)\times r$ submatrix of $\widetilde{B}$, $\be_0 = (e_1, \ldots, e_r)$, $\bx_N = (x_1, \ldots, x_r)$, $\bx_F = (x_{r+1}, \ldots, x_n )$, $\ind{\Gamma}^N R = \sum_{i=1}^{r}g_i[\Gamma_i]$ and $\ind{\Gamma}^F R = \sum_{i=r+1}^{n}g_i[\Gamma_i]$.  

It suffices to show that, if $\be$ is such that $\Gr{\be} \big( \Hom{\cC}(\Gamma, \Sigma R) \big)$ is non-empty, then  $\ind{\Gamma}^N R + B \cdot \be_0$ (the power of $\bx_N$) has at least one negative entry.

Assume first that $\be$ is non-zero.  Since $B$ is skew-symmetric, we have that $\be_0 \cdot B \cdot \be_0$ vanishes, and thus $\be_0 \cdot (\ind{\Gamma}^N R + B \cdot \be_0) = \be_0 \cdot \ind{\Gamma}^N R$.  Thus it is sufficient to prove that $\be_0 \cdot \ind{\Gamma}^N R$ is negative.  Since the last $(n-r)$ entries of $\be$ are zero, we have that $\be_0 \cdot \ind{\Gamma}^N R = \be \cdot \ind{\Gamma} R$.

Let $F$ be the functor $\Hom{\cC}(\Gamma, \Sigma (?))$.  Let $L$ be a submodule of $FR$ with dimension vector $\be$, and let $V$ be an object of $\cD$ such that $FV = L$.  Then, by \cite[Proposition 4.15]{Plamondon10} (or rather, by the second-to-last equality at the end of the proof)\footnote{There is a typo in this equality in \cite{Plamondon10}; the term $\ind{\Gamma}\Sigma Y$ should be read $\ind{\Gamma} Y$. }, we have that
\begin{eqnarray*}
  \dim (\Sigma \Gamma) (V, \Sigma R) &=& \dim \Hom{J(Q,W)}(FV, FR) + \be \cdot \ind{\Gamma} R.
\end{eqnarray*}
Since $FV = L$ is a submodule of $FR$ and since the latter is non-zero, we have that $\dim \Hom{J(Q,W)}(FV, FR)$ is positive.  Moreover, let $i:V\rightarrow R$ be a lift in $\cC_{Q,W}$ of the inclusion of $L$ in $FR$.  Embed $i$ in a triangle to get a diagram
\begin{displaymath}
	\xymatrix{ U\ar[r]^u & V \ar[r]^i\ar[d]^f & R\ar[r]\ar@{.>}[dl]^0 & \Sigma U \\
	                    &  \Sigma R          
	}
\end{displaymath}
where the upper row is a triangle, and where $u$ is in the ideal $(\Gamma)$, since $Fu = 0$ ($Fi$ being injective).  Now, if $f$ is in $(\Sigma \Gamma) (V, \Sigma R)$, then the composition $fu$ vanishes since $\Gamma$ is rigid, and thus $f$ factors through $R$.  But $R$ is also rigid; thus $f$ is zero.  Therefore $\dim (\Sigma \Gamma) (V, \Sigma R) = 0$, and $\be \cdot \ind{\Gamma} R$ is negative; this is what we wanted to prove.

Now, assume that $\be$ is zero.  Then we need to show that $\ind{\Gamma}^N R$ has a negative entry.  Let $R = T_1^{a_1} \oplus \ldots \oplus T_s^{a_s}$, with each $T_i$ indecomposable. By hypothesis, one of the $T_i$, with $i\in \{1, \ldots, r\}$, is not a direct factor of $\Gamma$; thus $\ind{\Gamma}^N T_i$ has a negative entry (because rigid objects whose indices have only non-negative entries are in $\add\Gamma$; this is a consequence of \cite[Proposition 3.1]{Plamondon10}).  Moreover, by \cite[Theorem 3.7(1)]{Plamondon10}, the indices of $T_1, \ldots, T_s$ are sign-coherent.  Thus $\ind{\Gamma}^N R = \sum_{i=1}^{s}a_i\ind{\Gamma}^N T_i$ has a negative entry.  

Another way of dealing with the case $\be = 0$ is by replacing $V$ by $R$ in the above equality, which becomes
\begin{eqnarray*}
  \dim (\Sigma \Gamma) (R, \Sigma R) &=& \dim \Hom{J(Q,W)}(FR, FR) + \dimv FR \cdot \ind{\Gamma} R.
\end{eqnarray*}
 Since $R$ is rigid, the left-hand side vanishes; since $R$ is not in $\add\Gamma$, $FR$ is not zero.  Thus we get
 \begin{displaymath}
	\dimv FR \cdot \ind{\Gamma} R = -\dim \Hom{J(Q,W)}(FR, FR) < 0,
\end{displaymath} 
and since the entries of $\dimv FR$ corresponding to frozen vertices vanish, this implies that $\ind{\Gamma}^N R$ has a negative entry. This finishes the proof.
}

\begin{corollary}\label{coro::linear independence}
If $B$ is any skew-symmetric matrix and $\widetilde{B}$ is an $n\times r$ integer matrix whose top $r\times r$ matrix is $B$, then the cluster algebra $\cA(\widetilde{B})$ has the proper Laurent monomial property.  In particular, the cluster monomials of any skew-symmetric cluster algebra of geometric type are linearly independent over $\bZ\bP$.
\end{corollary}
\demo{If the first statement of the theorem is true, then the independence of cluster monomials follows from \cite[Theorem 6.4]{CL11}.  Let us prove that $\cA(\widetilde{B})$ has the proper Laurent monomial property.

Let $Q$ be the quiver associated to $\widetilde{B}$, and let $W$ be a non-degenerate potential on $Q$.  Consider the cluster category $\cC_{Q,W}$.  Let $(u_1, \ldots, u_r)$ be a cluster of $\cA(\widetilde{B})$, and let $z = u_1^{a_1}\cdots u_r^{a_r}$ be a cluster monomial.  By \cite[Theorem 4.1]{Plamondon09}, there exists a rigid object $T = T_1\oplus \ldots \oplus T_n$ of $\cC_{Q,W}$ obtained from $\Gamma$ by iterated mutations at vertices in $\{1, \ldots, r \}$ and such that 
\begin{displaymath}
	z = CC(T_1^{a_1} \oplus \ldots \oplus T_r^{a_r}).
\end{displaymath}
Assume, moreover, that there is an index $i$ such that $u_i$ is not in the initial cluster and such that $a_i$ is positive.  Then $R = T_1^{a_1} \oplus \ldots \oplus T_r^{a_r}$ satisfies the hypothesis of Theorem \ref{theo::rigid}; thus $z$ is a $\bZ\bP$-linear combination of proper Laurent monomials in the initial cluster.  This finishes the proof.
}

\begin{corollary}\label{coro::bijection}
The surjection of Theorem \ref{theo::surjection} is a bijection.
\end{corollary}
\demo{Assume that $R$ and $R'$ are indecomposable reachable objects such that $CC(R) = CC(R')$.  Let $\bi$ be a sequence of vertices of $Q$ and $j$ be a vertex of $Q$ such that $\mu_{\bi}(R)$ is isomorphic to $\Gamma_j$.  Then we have that
\begin{displaymath}
	x_j = CC(\Gamma_j) = CC(\mu_{\bi}(R)) = \mu_{\bi}(CC(R)) = \mu_{\bi}(CC(R')) = CC(\mu_{\bi}(R')).
\end{displaymath}
Therefore, $CC(\mu_{\bi}(R'))$ is not a linear combination of proper Laurent monomials in the initial variables.  By Theorem \ref{theo::rigid}, this implies that $\mu_{\bi}(R')$ lies in $\add\Gamma$.  Since $CC(\mu_{\bi}(R')) = x_j$, we must have that $\mu_{\bi}(R')$ is isomorphic to $\Gamma_j$, and thus to   $\mu_{\bi}(R)$.  Therefore, $R$ and $R'$ are isomorphic.
}

As noticed in \cite[Theorem 4.1]{BMRT07} for acyclic cluster algebras with trivial coefficients, the bijection of Corollary \ref{coro::bijection} allows us to prove that seeds are determined by their clusters, as conjectured in \cite[Conjecture 4.14(2)]{FominZelevinsky03a}.  This was also proved in \cite{GSV08} for cluster algebras of geometric type, and for cluster algebras whose matrix is of full rank (including the skew-symmetrizable case).

\begin{corollary}
In any skew-symmetric cluster algebra, each seed is determined by its cluster.
\end{corollary}
\demo{ We first deal with cluster algebras with trivial coefficients.  Let $(\bu, Q)$ and $(\bu, Q')$ be two seeds with the same cluster.  By Corollary~\ref{coro::bijection}, there is a unique reachable object $T$ corresponding to $\bu$ in the cluster category; moreover, the ordinary quiver of $\End{\cC}(T)$ has to be both $Q$ and $Q'$.  Thus the two seeds coincide.  The result is proved for trivial coefficients.

If the coefficients are taken in an arbitrary semifield $\bP$, then the morphism of semifields $\bP \rightarrow \{1\}$ induces a morphism of cluster algebras $f:\cA \rightarrow \cA_{triv}$, where $\cA_{triv}$ is the cluster algebra defined from the same quiver as $\cA$, but with trivial coefficients.  Thus if two seeds $(\bu, Q)$ and $(\bu, Q')$ in $\cA$ have the same cluster, their images in $\cA_{triv}$ by $f$ also have the same cluster, and thus $Q$ and $Q'$ are equal, since clusters determine their seeds when the coefficients are trivial. 
}

To close this section, we present a proof of a maximality property that is satisfied by clusters, which was formulated as a question by Peter Tingley.

\begin{lemma}\label{lemm::invertible}
 Let $\bu$ be a cluster in a skew-symmetric cluster algebra $\cA$, and let $y$ be a cluster variable that is not in $\bu$.  Then $y$ is not invertible in $\bZ[u_1^{\pm 1}, \ldots, u_n^{\pm 1}]$ (in other words, $y$ is not a Laurent monomial in $\bu$).
\end{lemma}
\demo{ Assume that $y$ is a Laurent monomial in $\bu$.  Write $y = \pm u_1^{a_1}\cdots u_n^{a_n}$.  By Corollary \ref{coro::linear independence}, $y$ is a proper Laurent monomial; thus we can assume that $a_1$ is negative.  Let $u_1'$ be the cluster variable obtained by mutating $\bu$ at vertex $1$.  Then we get that
\begin{displaymath}
 y = \frac{(u_1')^{-a_1} u_2^{a_2} \cdots u_n^{a_n}}{(\prod_{1\rightarrow i} u_i + \prod_{j \rightarrow 1} u_j)^{-a_1} }.
\end{displaymath}
If the vertex $1$ is not incident with any arrow, then the denominator of this expression is a positive power of $2$, contradicting the fact that $y$ lies in $\bZ[u_1'^{\pm 1}, u_2^{\pm 1} \ldots, u_n^{\pm 1}]$.  Otherwise the denominator is a sum of two distinct monomials, contradicting the Laurent phenomenon.  Thus $y$ cannot be invertible.
}

\begin{proposition}\label{prop::invertible}
 Let $\{z_1, \ldots, z_p\}$ be a set of cluster variables in a skew-symmetric cluster algebra.  Assume that arbitrary finite products of elements of this set are cluster monomials, and that the set is maximal for this property.  Then the set is a cluster.
\end{proposition}
\demo{ By our assumption, the product $z_1\cdots z_p$ is a cluster monomial; thus there exists a cluster $\bu$ and non-negative integers $a_1, \ldots, a_n$ such that
\begin{displaymath}
 z_1\cdots z_p = u_1^{a_1} \cdots u_n^{a_n}.
\end{displaymath}
This implies that $z_1, \ldots, z_p$ are all invertible in $\bZ[u_1^{\pm 1}, \ldots, u_n^{\pm 1}]$; by Lemma \ref{lemm::invertible}, they must all belong to $\bu$.  Therefore, the set $\{z_1, \ldots, z_p\}$ is a subset of the cluster $\bu$.  By maximality, the two sets are equal.
}

\section{Exchange graphs}

\begin{definition}[Definition 4.2 of \cite{FZ07}]
The \emph{exchange graph} of a cluster algebra $\cA(\bx, \by, B)$ is the $n$-regular connected graph whose vertices are the seeds of $\cA(\bx, \by, B)$ (up to simultaneous renumbering of rows, columns and
variables) and whose edges connect the seeds related by a single mutation.
\end{definition}

\begin{conjecture}[Conjecture 4.3 of \cite{FZ07}]\label{conj::exchange}
The exchange graph of a cluster algebra $\cA(\bx, \by, B)$ only depends on the matrix $B$.
\end{conjecture}

\begin{definition}
Let $(Q,W)$ be a non-degenerate quiver with potential (where the vertices $r+1, \ldots, n$ are frozen), and let $\cC_{Q,W}$ be the associated cluster category.
\begin{enumerate}
	\item The \emph{exchange graph of $\Gamma$} is the $r$-ary graph whose vertices are isoclasses of objects reachable from $\Gamma$ (by mutation at non-frozen vertices), and where two vertices are joined by an edge if the corresponding objects are related by a mutation.
	
	\item The \emph{exchange graph of indices} is the $r$-ary graph obtained from the exchange graph of $\Gamma$ by replacing any vertex $T = T_1\oplus \ldots \oplus T_n$ by the $n$-tuple $(\ind{\Gamma}T_1, \ldots, \ind{\Gamma}T_n)$ of elements of $K_0(\add\Gamma) \cong \bZ^n$.
	
	\item The \emph{exchange graph of reduced indices} is the graph obtained from the exchange graph of indices by replacing each $n$-tuple 
	\[
	(\ind{\Gamma}T_1, \ldots, \ind{\Gamma}T_n)
	\]
	 by the $r$-tuple $(\ind{\Gamma}^N T_1, \ldots, \ind{\Gamma}^N T_r)$, where $\ind{\Gamma}^N T_i$ is obtained from $\ind{\Gamma}T_i$ by forgetting the coefficients of $[\Gamma_{r+1}], \ldots, [\Gamma_{n}]$.
\end{enumerate}
\end{definition}

The following corollary is a consequence of Corollary \ref{coro::bijection}.  
\begin{corollary}\label{coro::exchange}
The exchange graphs of $\cA(\widetilde{B})$, of $\Gamma$, of indices and of reduced indices are isomorphic.
\end{corollary}
\demo{ The fact that the exchange graphs of $\cA(\widetilde{B})$ and of $\Gamma$ are isomorphic is a direct consequence of Corollary \ref{coro::bijection}.  The fact that the exchange graphs of $\Gamma$ and of indices are isomorphic follows from the fact that rigid objects are determined by their index, by section 3.1 of \cite{Plamondon10}, \emph{cf.} also \cite{DK08}.

Let $\Gamma_N$ be the direct sum of the $\Gamma_j$, $1\leq j \leq r$, and $\Gamma_F$ be the direct sum of the $\Gamma_j$, $r+1\leq j \leq n$.
To prove that the exchange graph of indices and of reduced indices are isomorphic, it is sufficient to show that if $X$ and $X'$ are indecomposable reachable objects (not in $\add\Gamma_F$) such that $\ind{\Gamma}^NX = \ind{\Gamma}^NX'$, then they are isomorphic.  Let
\begin{displaymath}
	T_1 \rightarrow T_0 \oplus R \rightarrow X \rightarrow \Sigma T_1 \textrm{  and  }  T_1 \rightarrow T_0 \oplus R' \rightarrow X' \rightarrow \Sigma T_1
\end{displaymath}
be triangles, with $T_1, T_0 \in \add\Gamma_N$ and $R,R' \in \add\Gamma_F$.  We can assume that $T_1$ has no direct factor in $\add\Gamma_F$ since we have that $\Hom{\cC}(X, \Sigma \Gamma_F)$ and $\Hom{\cC}(X', \Sigma \Gamma_F)$ vanish (this is true because we only allow mutations at non-frozen vertices).  But then $X\oplus R'$ and $X'\oplus R$ are rigid objects which have the same index; by \cite{DK08}, they are isomorphic.  Since $X$ and $X'$ are indecomposable and not in $\add\Gamma_F$, we get that they are isomorphic.
}

\begin{lemma}\label{lemm::reduced indices}
Let $(Q,W)$ be a non-degenerate quiver with potential, and let $Q^0$ be the quiver obtained from $Q$ by deleting the frozen vertices and arrows incident with them.  Then the graph of reduced indices of $(Q,W)$ only depends on $Q^0$.
\end{lemma}
\demo{ Note first that it follows from the definition of mutation of quivers that $\mu_{i}(Q)^0 = \mu_{i}(Q^0)$ for any non-frozen vertex $i$.

For any quiver $R$ and any sequence $\bi = (i_1, \ldots, i_s)$ of non-frozen vertices, denote by $\mu_{\bi}(R)$ the mutated quiver $\mu_{i_s}\cdots\mu_{i_2}\mu_{i_1}(R)$.  

Let $\bi = (i_1, \ldots, i_s)$ be a sequence of non-frozen vertices of $Q$, and let $T = T_1\oplus \ldots \oplus T_n$ be the object obtained from $\Gamma$ by iterated mutations at $\bi$.  The corresponding vertex in the exchange graph of reduced indices is $(\ind{\Gamma}^N T_1, \ldots, \ind{\Gamma}^N T_r)$.  We must prove that each entry only depends on $Q^0$ (and, implicitly, on the sequence $\bi$ of mutations).

Recall from \cite{KY09} the functor $\Phi_i : \cC_{\mu_i(Q,W)} \rightarrow \cC_{Q,W}$ which sends the direct factor $(\Gamma_{\mu_i(Q,W)})_j$ of $\Gamma_{\mu_i(Q,W)}$ to $\Gamma_j$ if $i\neq j$ and to a cone of the natural morphism
\begin{displaymath}
	\Gamma_i \longrightarrow \bigoplus_{\alpha:i\rightarrow k}\Gamma_k
\end{displaymath}
if $i = j$.    Following \cite{KY09}, $T$ is the image of $\Gamma_{\mu_{\bi}(Q,W)} \in \cC_{\mu_{\bi}(Q,W)}$ by the composition of functors
\begin{displaymath}
	\xymatrix{\cC_{\mu_{\bi}(Q,W)} \ar[r]^{\phantom{xx}\Phi_{i_s}} & \ldots \ar[r]^{\Phi_{i_2}\phantom{xx}} &  \cC_{\mu_{i_1}(Q,W)} \ar[r]^{\phantom{x}\Phi_{i_1}} & \cC_{Q,W}.
	}
\end{displaymath}
We will write $\Gamma' = \Gamma_{\mu_i(Q,W)}$ and $T = \Phi_{\bi}(\Gamma')$.

The lemma will be a consequence of the following stronger statement: \emph{If $X$ is a rigid object of $\cC_{\mu_\bi(Q,W)}$ which admits an index, then $\ind{\Gamma}^N \Phi_{\bi}(X)$ only depends on $Q^0$ and on $\ind{\Gamma'}^N X$}.

We prove this statement by induction on $s$.  For $s=1$, we only perform one mutation at vertex $i=i_1$.  
Let $\ind{\Gamma}\Phi_i(X) = \sum_{j=1}^n g_j[\Gamma_j]$ and $\ind{\Gamma'}X = \sum_{j=1}^n g'_j [\Gamma'_j]$.  Then the proof of \cite[Conjecture 7.12]{FZ07} given in \cite[Theorem 3.7 (iv)]{Plamondon10} gives us the formula
\begin{equation*}
g_j =
\begin{cases}
-g'_i & \text{if } i=j,\\
g'_j + g'_i [b'_{ji}]_+ - b'_{ji} \min(g'_i, 0) & \text{if } i\neq j,
\end{cases}
\end{equation*}
where $B' = (b'_{ij})$ is the matrix of $\mu_i(Q)$, and $[z]_+$ stands for $\max(z, 0)$ for any real number $z$.  It follows from this formula that $\ind{\Gamma}^N \Phi_i(X) = \sum_{j=1}^r g_j[\Gamma_j]$ only depends on $\ind{\Gamma'}^N X$ and on $\mu_i(Q)^0$, and this last quiver only depends on $Q^0$.  The statement is shown for $s=1$.

Assume it is proved for $s-1$. Let $\bi = (i_1, \ldots, i_s)$. Then $\ind{\Gamma}^N \Phi_{\bi}(X)$ only depends on $Q^0$ and on $\ind{\Gamma''}^N \Phi_{i_s}(X)$, where $\Gamma'' = \Gamma_{\mu_{i_{s-1}}\cdots\mu_{i_1}(Q,W)}$.  Moreover, $\ind{\Gamma''}^N \Phi_{i_s}(X)$ only depends on $\ind{\Gamma'}^N X$ and on $\mu_{i_{s-1}}\cdots\mu_{i_1}(Q)^0$.  This last quiver only depends on $Q^0$.   This proves the statement.

To finish the proof of the lemma, apply the statement to $X = \Gamma'_j$.  Since $\Phi_{\bi}(\Gamma'_j) = T_j$, and since $\Gamma'_j$ is determined by its index, we get that $\ind{\Gamma}^N T_j$ only depends on $Q_0$.  This finishes the proof of the lemma.
}

\begin{theorem}\label{theo::exchange graph}
Conjecture \ref{conj::exchange} is true for skew-symmetric cluster algebras of geometric type, that is, the exchange graph of $\cA(\widetilde{B})$ only depends on the top $r\times r$ matrix $B$.
\end{theorem}
\demo{ By Corollary \ref{coro::exchange}, the exchange graph of $\cA(\widetilde{B})$ is isomorphic to the exchange graph of reduced indices, which only depends on $B$ by Lemma \ref{lemm::reduced indices}.
}

Using similar arguments, we can prove a stronger result concerning the cluster complex.

\begin{definition}[\cite{FZ03}]
 The \emph{cluster complex} of a cluster algebra $\cA(\bx, \by, B)$ is the simplicial complex whose set of vertices is the set of cluster variables and whose simplices are the subsets of clusters.
\end{definition}

\begin{theorem}\label{theo::cluster complexes}
 Let $\cA(\widetilde{B})$ be a cluster algebra of geometric type.  Then the cluster complex of $\cA(\widetilde{B})$ only depends on the top $r\times r$ matrix $B$. 
\end{theorem}
\demo{ The proof of this result follows the lines of the proof of Theorem \ref{theo::exchange graph}.  Let $Q$ be the quiver (with frozen vertices) associated to $\widetilde{B}$, and let $W$ be a non-degenerate potential on $Q$. Define the following simplicial complexes:
\begin{enumerate}
 \item the \emph{reachable complex} of $\cC_{Q,W}$, whose vertices are indecomposable reachable objects and whose simplices correspond to direct factors of reachable cluster-tilting objects;
 \item the \emph{index complex}, obtained from the reachable complex by replacing each vertex by the index of the corresponding object;
 \item the \emph{reduced index complex}, obtained from the index complex by replacing each vertex of the form $\ind{\Gamma}X$ by $\ind{\Gamma}^N X$.
\end{enumerate}
It follows from Corollary \ref{coro::bijection} that the cluster complex of $\cA(\widetilde{B})$ and the reachable complex of $\cC_{Q,W}$ are isomorphic.  Since rigid objects are determined by their index, the reachable complex is isomorphic to the index complex.  Moreover, since a reachable indecomposable object $X$ is in fact determined by $\ind{\Gamma}^N X$ (see the proof of Corollary \ref{coro::exchange}), the index complex is isomorphic to the reduced index complex.

Notice that the reduced index complex is completely determined by the maximal simplex corresponding to the initial cluster, and by the mutation rule.  Indeed, from there one can construct the rest of the simplicial complex by iterated mutations.  Moreover, it was proved in Lemma \ref{lemm::reduced indices} that the mutation rule for reduced indices only depends on the top $r\times r$ matrix $B$.  This finishes the proof.
}

\section{Cluster algebras with arbitrary coefficients}
In this section, we extend Corollary \ref{coro::linear independence} and Theorem \ref{theo::exchange graph} to the case of skew-symmetric cluster algebras with arbitrary coefficients.  For the definition of cluster algebras with arbitrary coefficients, we refer the reader to \cite[Section 2]{FZ07}.

The linear independence of cluster monomials will follow from Corollary \ref{coro::linear independence} and from the ``separation of additions'' of Fomin-Zelevinsky, which we recall.  

\begin{theorem}[Theorem 3.7 of \cite{FZ07}]\label{theo::separation additions}
Let $\cA$ be a cluster algebra associated to the matrix $B$ and with coefficients in an arbitrary semifield $\bP$.  Let $t$ be a seed of $\cA$, and let $x_{\ell, t}$ be a cluster variable of that seed.  Then
\begin{displaymath}
	x_{\ell, t} = \frac{X_{\ell, t}|_{\cF}(x_1, \ldots, x_m; y_1, \ldots, y_m)}{F_{\ell,t}|_{\bP}(y_1, \ldots, y_m)}.
\end{displaymath}
\end{theorem}
In this theorem, $X_{\ell, t}$ and $F_{\ell, t}$ are rational functions which depend on the cluster algebra with principal coefficients (see \cite[Definition 3.1]{FZ07}) and with defining matrix $B$.  They are defined in \cite[Definition 3.3]{FZ07}.

\begin{corollary}\label{coro::principal to arbitrary}
If the cluster algebra with principal coefficients $\cA_0$ associated to the matrix $B$ has the proper Laurent monomial property, then any cluster algebra associated to $B$ has the same property.
\end{corollary}
\demo{By definition, $X_{\ell, t}$ is an expression of a cluster variable with respect to the initial variables, inside the cluster algebra with principal coefficients.  By assumption, this is an initial variable or a proper Laurent monomial in the initial variables.  Specializing the coefficients to $y_1, \ldots, y_m$ and dividing by $F_{\ell,t}|_{\bP}(y_1, \ldots, y_m)$ does not change this property.  Therefore any cluster variable $x_{\ell, t}$ is a proper Laurent monomial.
}

\begin{corollary}
A skew-symmetric cluster algebra with arbitrary coefficients has the proper Laurent monomial property.  In particular, its cluster monomials are linearly independent.
\end{corollary}
\demo{ By Corollary \ref{coro::linear independence}, any skew-symmetric cluster algebra with principal coefficients has the proper Laurent monomial property (since 
cluster algebras with principal coefficients are of geometric type).  Therefore, by Corollary \ref{coro::principal to arbitrary}, any skew-symmetric cluster algebra also has this property.  Then Theorem~6.4 of \cite{CL11} implies that cluster monomials 
in such a cluster algebra are linearly independent.  Note that 
Theorem~6.4 of \cite{CL11} is only stated for cluster algebras of geometric type; however, its proof adapts without change to arbitrary cluster algebras.
}

The proof that the exchange graph of a cluster algebra only depends on the defining matrix will follow from Theorem \ref{theo::exchange graph} and from the following theorem of Fomin-Zelevinsky.

\begin{theorem}[Theorem 4.6 of \cite{FZ07}]\label{theo::covering}
The exchange graph of an arbitrary cluster algebra is covered by the exchange graph of the cluster algebra with principal coefficients having the same defining matrix.
\end{theorem}

\begin{corollary}
The exchange graph of any skew-symmetric cluster algebra only depends on its defining matrix.
\end{corollary}
\demo{ Let $\cA$ be any skew-symmetric cluster algebra, with defining matrix $B$.  Let $\cA_0$ and $\cA_{triv}$ be the cluster algebras with principal and trivial coefficients, respectively, associated to $B$.  Let $G$, $G_0$ and $G_{triv}$ be the exchange graphs of $\cA$, $\cA_0$ and $\cA_{triv}$, respectively.  By Theorem \ref{theo::covering}, $G_0$ is a covering of $G$.  Moreover, $G$ is a covering of $G_{triv}$.  Since $\cA_0$ and $\cA_{triv}$ are both of geometric type, Theorem \ref{theo::exchange graph} implies that the the composed covering
$G_0 \to G_{triv}$ is an isomorphism.  
Therefore, the coverings $G_0 \to G$ and $G \to G_{triv}$ are
isomorphisms.  This finishes the proof.
}

Similarly, we can prove that the cluster complex only depends on the defining matrix.

\begin{corollary}
 The cluster complex of any skew-symmetric cluster algebra only depends on its defining matrix.
\end{corollary}
\demo{Let $\cA$ be any skew-symmetric cluster algebra, with defining matrix $B$.  Let $\cA_0$ and $\cA_{triv}$ be the cluster algebras with principal and trivial coefficients, respectively, associated to $B$. Let $\Delta$, $\Delta_0$ and $\Delta_{triv}$ be the cluster complexes of $\cA$, $\cA_0$ and $\cA_{triv}$, respectively.

First, notice that the cluster complex of $\cA_0$ covers the cluster complexes of both $\cA$ and $\cA_{triv}$.  This follows from the separation of additions \cite[Theorem 3.7]{FZ07}(the proof of this statement is in fact analogous to the proof of \cite[Theorem 4.6]{FZ07}).  Moreover, the cluster complex of $\cA$ covers that of $\cA_{triv}$.  Since $\cA_0$ and $\cA_{triv}$ are both of geometric type, Theorem \ref{theo::cluster complexes} implies that the composed covering $\Delta_0 \to \Delta_{triv}$ is an isomorphism.  Therefore, the coverings $\Delta_0 \to \Delta$ and $\Delta \to \Delta_{triv}$ are
isomorphisms.  This finishes the proof.
}

\section*{Acknowledgements}
The authors would like to thank Peter Tingley for his question which led to Proposition \ref{prop::invertible}, and David Hernandez for informing them about this question.

\bibliographystyle{amsplain} 

\providecommand{\bysame}{\leavevmode\hbox to3em{\hrulefill}\thinspace}
\providecommand{\MR}{\relax\ifhmode\unskip\space\fi MR }
\providecommand{\MRhref}[2]{%
  \href{http://www.ams.org/mathscinet-getitem?mr=#1}{#2}
}
\providecommand{\href}[2]{#2}

\end{document}